\documentclass[a4paper,12pt]{article}

\usepackage{amsfonts}
\usepackage{amssymb}
\usepackage{latexsym}
\usepackage{amsmath}
\usepackage{amsthm}
\usepackage{graphicx}
\usepackage{color}

\newtheorem{theorem}{Theorem}

\newtheorem{definition}{Definition}
\newtheorem{remark}{Remark}

\newtheorem{conjecture}{Conjecture}

\newcommand{\NN}{\mathbb{N}}
\newcommand{\N}{\mathbb{N}}
\newcommand{\RR}{\mathbb{R}}
\newcommand{\ZZ}{\mathbb{Z}}
\newcommand{\bsgamma}{\boldsymbol{\gamma}}
\newcommand{\bsalpha}{\boldsymbol{\alpha}}
\newcommand{\bsx}{\boldsymbol{x}}
\newcommand{\bsy}{\boldsymbol{y}}
\newcommand{\bsf}{\boldsymbol{f}}
\newcommand{\bsj}{\boldsymbol{j}}
\newcommand{\bsm}{\boldsymbol{m}}
\newcommand{\bst}{\boldsymbol{t}}
\newcommand{\bsone}{\boldsymbol{1}}
\newcommand{\bszero}{\boldsymbol{0}}

\newcommand{\cS}{\mathcal{S}}

\newcommand{\rd}{\,\mathrm{d}}

\newcommand{\FF}{\mathbb{F}}

\begin{document}

\title{Discrepancy of Digital Sequences: New Results on a Classical QMC Topic}
\author{Friedrich Pillichshammer\thanks{F. Pillichshammer is supported by the Austrian Science Fund (FWF) Project F5509-N26, which is a part of the Special Research Program ``Quasi-Monte Carlo Methods: Theory and Applications''.}}

\date{}

\maketitle

\begin{center}
{\bf Dedicated to Henri Faure on the occasion of his $80^{{\rm th}}$ birthday.}
\end{center}
\vspace{1cm}

\begin{abstract}
The theory of digital sequences is a fundamental topic in QMC theory. Digital sequences are prototypes of sequences with low discrepancy. First examples were given by Il'ya Meerovich Sobol' and by Henri Faure with their famous constructions. The unifying theory was developed later by Harald Niederreiter. Nowadays there is a magnitude of examples of digital sequences and it is classical knowledge that the star discrepancy of the initial $N$ elements of such sequences can achieve a rate of order $(\log N)^s/N$, where $s$ denotes the dimension. On the other hand, very little has been known about the $L_p$ norm of the discrepancy function of digital sequences for finite $p$, apart from evident estimates in terms of star discrepancy.
In this article we give a review of some recent results on various types of discrepancy of digital sequences. This comprises: star discrepancy and weighted star discrepancy, $L_p$-discrepancy, discrepancy with respect to bounded mean oscillation and exponential Orlicz norms, as well as Sobolev, Besov and Triebel-Lizorkin norms with dominating mixed smoothness. 
\end{abstract}

\section*{Preamble}

This paper is devoted to Henri Faure who celebrated his $80^{{\rm th}}$ birthday on July 12, 2018. Henri is well known for his pioneering work on low-discrepancy sequences. As an example we would like to mention his famous paper \cite{fau} from 1982 in which he gave one of the first explicit constructions of digital sequences in arbitrary dimension with low star discrepancy. These sequences are nowadays known as {\it Faure sequences}. 

I met Henri for the first time at the MCQMC conference 2002 in Singapore. Later, during several visits of Henri in Linz, we started a fruitful cooperation which continues to this day. I would like to thank Henri for this close cooperation and for his great friendship and wish him and his family all the best for the future. 

\section{Introduction}
\label{sec:1}

We consider infinite sequences $\mathcal{S}=(\bsx_n)_{n \ge 0}$ of points $\bsx_n$ in the $s$-dimensional unit cube $[0,1)^s$. For $N \in \mathbb{N}$ let $\mathcal{S}_N=(\bsx_n)_{n=0}^{N-1}$ be the initial segment of $\mathcal{S}$ consisting of the first $N$ elements. 

According to Weyl~\cite{weyl} a sequence $\mathcal{S}=(\bsx_n)_{n \ge 0}$ is uniformly distributed (u.d.) if for every axes-parallel box $J \subseteq [0,1)^s$ it is true that $$\lim_{N \rightarrow \infty} \frac{\#\{n \in \{0,\ldots,N-1\} \ : \ \bsx_n \in J\}}{N}={\rm Volume}(J).$$ An extensive introduction to the theory of uniform distribution of sequences can be found in the book of Kuipers and Niederreiter~\cite{kuinie}.

There are several equivalent definitions of uniform distribution of a sequence and one of them is of particular importance for quasi-Monte Carlo (QMC) integration. Weyl proved that a sequence $\mathcal{S}$ is u.d. if and only if for every Riemann-integrable function $f:[0,1]^s \rightarrow \RR$ we have 
\begin{equation}\label{udfa}
\lim_{N \rightarrow \infty} \frac{1}{N}\sum_{n=0}^{N-1} f(\bsx_n)=\int_{[0,1]^s} f(\bsx) \, {\rm d}\bsx.
\end{equation}
The average of function evaluations on the left-hand side is nowadays called a {\it QMC rule}, $$Q_N(f)=\frac{1}{N}\sum_{n=0}^{N-1} f(\bsx_n).$$ Hence, in order to have a QMC rule converging to the true value of the integral of a function it has to be based on a u.d. sequence. A quantitative version of \eqref{udfa} can be stated in terms of discrepancy.

\begin{definition}
For a finite initial segment $\mathcal{S}_N$ of a sequence (or a finite point set) in $[0,1)^s$ the {\it local discrepancy function} $\Delta_{\mathcal{S}_N}:[0,1]^s \rightarrow \RR$ is defined as $$\Delta_{\mathcal{S}_N}(\bst)=\frac{\#\{n \in \{0,1,\ldots,N-1\}\ : \ \bsx_n \in [\bszero,\bst)\}}{N}- t_1 t_2 \cdots t_s,$$ where $\bst=(t_1,t_2,\ldots,t_s)$, $[\bszero,\bst)=[0,t_1)\times [0,t_2)\times \ldots \times [0,t_s)$, and hence $t_1 t_2\cdots t_s={\rm Volume}([\bszero,\bst))$.

For $p \ge 1$ the {\it $L_p$ discrepancy} of $\mathcal{S}_N$ is defined as the $L_p$ norm of the local discrepancy function $$L_{p,N}(\mathcal{S}_N)=\|\Delta_{\mathcal{S}_N}\|_{L_p([0,1]^s)}=\left(\int_{[0,1]^s} |\Delta_{\cS_N}(\bst)|^p \rd \bst\right)^{1/p}$$ with the usual adaptions if $p=\infty$. In this latter case one often talks about {\it star discrepancy} which is denoted by $D_N^*(\mathcal{S}_N):=L_{\infty,N}(\mathcal{S}_N)$.

For an infinite sequence $\mathcal{S}$ in $[0,1)^s$ we denote the $L_p$ discrepancy of the first $N$ points by $L_{p,N}(\mathcal{S})= L_{p,N}(\mathcal{S}_N)$ for $N \ge 1$. 
\end{definition}

It is well-known that a sequence $\mathcal{S}$ is u.d. if and only if $\lim_{N \rightarrow \infty} L_{p,N}(\mathcal{S})=0$ for some $p \ge 1$. A quantitative version of \eqref{udfa} is the famous {\it Koksma-Hlawka inequality} which states that for every function $f:[0,1]^s \rightarrow \RR$ with bounded variation $V(f)$ in the sense of Hardy and Krause and for every finite sequence $\mathcal{S}_N$ of points in $[0,1)^s$ we have $$\left|\int_{[0,1]^s} f(\bsx) \, {\rm d}\bsx- \frac{1}{N}\sum_{n=0}^{N-1} f(\bsx_n)\right| \le V(f) D_N^*(\mathcal{S}_N).$$ The Koksma-Hlawka inequality is the fundamental error estimate for QMC rules and the basis for QMC theory.  Nowadays there exist several versions of this inequality which may also be based on the $L_p$ discrepancy or other norms of the local discrepancy function. One often speaks about ``Koksma-Hlawka type inequalities''. For more information and for introductions to QMC theory we refer to \cite{DKS,DP2010,LP14_book,niesiam}. 

It is clear that QMC requires sequences with low discrepancy in some sense and this motivates the study of ``low discrepancy sequences''. On the other hand discrepancy is also an interesting topic by itself that is intensively studied (see, e.g., the books \cite{BC,CST,DT,DP2010,kuinie,Mat99,niesiam}). 

In the following we collect some well-known facts about $L_p$ discrepancy of finite and infinite sequences.

\section{Known facts about the $L_p$ discrepancy}

We begin with results on finite sequences: for every $p \in(1,\infty]$  and $s \in \NN$ there exists a $c_{p,s}>0$ such that for every finite $N$-element sequence $\cS_N$ in $[0,1)^s$ with $N \ge 2$ we have 
\begin{equation*}
L_{p,N}(\cS_N) \ge c_{p,s} \frac{(\log N)^{\frac{s-1}{2}}}{N} \ \ \mbox{ and } \ \  D_N^{\ast}(\cS_N) \ge c_{\infty,s} \frac{(\log N)^{\frac{s-1}{2}+\eta_s}}{N}
\end{equation*}
for some $\eta_s \in (0,\tfrac{1}{2})$. The result on the left hand side for $p \ge 2$ is a celebrated result by Roth~\cite{R54} from 1954 that was extended later by Schmidt~\cite{S77} to the case $p \in (1,2)$. The general lower bound for the star discrepancy is an important result of Bilyk, Lacey and Vagharshakyan~\cite{BLV08} from 2008. As shown by Hal\'{a}sz~\cite{H81}, the $L_p$ estimate is also true for $p=1$ and $s=2$, i.e.,  there exists a positive constant $c_{1,2}$ with the following property:  for every finite sequence $\cS_N$ in $[0,1)^2$ with $N \ge 2$ we have 
\begin{equation}\label{lbdl1D2dipts}
L_{1,N}(\cS_N) \ge c_{1,2} \frac{\sqrt{\log N}}{N}.
\end{equation} 
Schmidt showed for $s=2$ the improved lower bound on star discrepancy $$D_N^*(\cS_N) \ge c_{\infty,2} \frac{\log N}{N}$$ for some $c_{\infty,2}>0$. On the other hand, it is known that for every $s,N \in \NN$ there exist finite sequences $\cS_N$ in $[0,1)^s$ such that 
$$D_N^{\ast}(\cS_N) \lesssim_s \frac{(\log N)^{s-1}}{N}.$$ First examples for such sequences are the Hammersley point sets, see, e.g., \cite[Section~3.4.2]{DP2010} or \cite[Section~3.2]{niesiam}. 

Similarly, for every $s,N \in \NN$ and every $p \in [1,\infty)$ there exist finite sequences $\cS_N$ in $[0,1)^s$ such that 
\begin{equation}\label{uplpps}
L_{p,N}(\cS_N) \lesssim_{s,p} \frac{(\log N)^{\frac{s-1}{2}}}{N}.
\end{equation}
Hence, for $p \in (1,\infty)$ and arbitrary $s \in \NN$ we have matching lower and upper bounds. For both $p=1$ and $p=\infty$ we have matching lower and upper bounds only for $s=2$. The result in \eqref{uplpps} was proved by Davenport \cite{D56} for $p= 2$, $s= 2$, by Roth \cite{R80} for $p= 2$ and arbitrary $s$ and finally by Chen \cite{C80} in the general case. Other proofs were found by Frolov~\cite{Frolov}, Chen~\cite{C83}, Dobrovol'ski{\u\i}~\cite{Do84}, Skriganov~\cite{Skr89, Skr94}, Hickernell and Yue~\cite{HY00}, and Dick and Pillichshammer~\cite{DP05b}. For more details on the early history of the subject see the monograph \cite{BC}. Apart from Davenport, who gave an explicit construction in dimension $s=2$, these results are pure existence results and explicit constructions of point sets were not known until the beginning of this millennium. First explicit constructions of point sets with optimal order of $L_2$ discrepancy in arbitrary dimensions have been provided in 2002 by Chen and Skriganov \cite{CS02} for $p=2$ and in 2006 by Skriganov \cite{S06} for general $p$. Other explicit constructions are due to Dick and Pillichshammer \cite{DP14a} for $p=2$, and Dick \cite{D14} and Markhasin \cite{M15} for 
general $p$.

Before we summarize results about infinite sequences some words about the conceptual difference between the discrepancy of finite and infinite sequences are appropriate. Matou\v{s}ek~\cite{Mat99} explained this in the following way: while for finite sequences one is interested in the distribution behavior of the whole sequence $(\bsx_0,\bsx_1,\ldots,\bsx_{N-1})$ with a fixed number of elements $N$, for infinite sequences one is interested in the discrepancy of all initial segments $(\bsx_0)$, $(\bsx_0,\bsx_1)$, $(\bsx_0,\bsx_1,\bsx_2)$, \ldots, $(\bsx_0,\bsx_1,\bsx_2,\ldots,\bsx_{N-1})$, simultaneously for $N \in \NN$. In this sense the discrepancy of finite sequences can be viewed as a static setting and the discrepancy of infinite sequences as a dynamic setting. 

Using a method from Pro{\u\i}nov~\cite{Pro86} (see also \cite{DP14a}) the results about lower bounds on $L_p$ discrepancy for finite sequences can be transferred to the following lower bounds for infinite sequences: for every $p \in(1,\infty]$ and every $s \in \NN$ there exists a $C_{p,s}>0$ such that for every infinite sequence $\cS$ in $[0,1)^s$  
\begin{equation}\label{lbdlpdiseq}
L_{p,N}(\cS) \ge C_{p,s} \frac{(\log N)^{\frac{s}{2}}}{N} \ \ \ \ \mbox{infinitely often}
\end{equation}
and 
\begin{equation}\label{bdstdisequ}
D_N^{\ast}(\cS) \ge C_{\infty,s} \frac{(\log N)^{\frac{s}{2}+\eta_s}}{N} \ \ \ \ \mbox{infinitely often,}
\end{equation}
where $\eta_s \in (0,\tfrac{1}{2})$ is independent of the concrete sequence. For $s=1$ the result holds also for the case $p=1$, i.e., for every $\cS$ in $[0,1)$ we have
\begin{equation*}
L_{1,N}(\cS) \ge c_{1,1} \frac{\sqrt{\log N}}{N} \ \ \ \ \mbox{infinitely often,}
\end{equation*}
and the result on the star discrepancy can be improved to (see Schmidt \cite{S72}; see also \cite{B82,Lar15,LarPu15})
\begin{equation}\label{lpseqschm}
D_N^{\ast}(\cS) \ge c_{\infty,1} \frac{\log N}{N} \ \ \ \ \mbox{infinitely often.}
\end{equation}

On the other hand, for every dimension $s$ there exist infinite sequences $\cS$ in $[0,1)^s$ such that 
\begin{equation}\label{bdlds}
D_N^{\ast}(\cS)\lesssim_s \frac{(\log N)^s}{N} \ \ \ \ \mbox{for all $N \ge 2$.}
\end{equation}
Informally one calls a sequence a {\it low-discrepancy sequence} if its star discrepancy satisfies the bound \eqref{bdlds}. Examples of low-discrepancy sequences are: 
\begin{itemize}
\item Kronecker sequences $(\{n \bsalpha\})_{n \ge 0}$, where $\bsalpha \in \RR^s$ and where the fractional part function $\{\cdot\}$ is applied component-wise. In dimension $s=1$ and if $\alpha \in \RR$ has bounded continued fraction coefficients, then the Kronecker sequence $(\{n \alpha\})_{n \ge 0}$ has star discrepancy of exact order of magnitude $\log N/N$; see \cite[Chapter~3]{niesiam} for more information.
\item Digital sequences: the prototype of a digital sequence is the van der Corput sequence in base $b$ which was introduced by van der Corput~\cite{vdc35} in 1935. For an integer $b \ge 2$ (the ``basis'') the $n^{{\rm th}}$ element of  this sequence is given by $x_n=n_0 b^{-1}+n_1 b^{-2}+ n_2 b^{-3}+\cdots$ whenever $n$ has $b$-adic expansion $n=n_0+n_1 b+n_2 b^2+\cdots$. The van der Corput sequence has star discrepancy of exact order of magnitude $\log N/N$; see the recent survey article \cite{FKP} and the references therein.  

Multi-dimensional extensions of the van der Corput sequence are the Halton sequence \cite{H60}, which is the component-wise concatenation of van der Corput sequences in pairwise co-prime bases, or digital $(t,s)$-sequences, where the basis $b$ is the same for all coordinate directions. First examples of such sequences have been given by Sobol' \cite{sob} and by Faure \cite{fau}. Later the general unifying concept has been introduced by Niederreiter~\cite{N87} in 1987. Halton sequences in pairwise co-prime bases as well as digital $(t,s)$-sequences have star discrepancy of order of magnitude of at most $(\log N)^s/N$; see Section~\ref{digtssequ}.  
\end{itemize}

Except for the one-dimensional case, there is a gap for the $\log N$ exponent in the lower and upper bound on the star discrepancy of infinite sequences (cf. Eq.~\eqref{bdstdisequ} and Eq.~\eqref{bdlds}) which seems to be very difficult to close. There is a grand conjecture in discrepancy theory which share many colleagues (but it must be mentioned that there are also other opinions; see, e.g., \cite{BL14}):
\begin{conjecture}\label{con1}
For every $s \in \NN$ there exists a $c_s>0$ with the following property: for every $\cS$ in $[0,1)^s$ it holds true that $$D_N^{\ast}(\cS) \ge c_s \frac{(\log N)^s}{N}\ \ \ \ \mbox{infinitely often.}$$ 
\end{conjecture}

For the $L_p$ discrepancy of infinite sequences with finite $p$ the situation is different. It was widely assumed that the general lower bound of Roth-Schmidt-Pro{\u\i}nov in Eq. \eqref{lbdlpdiseq} is optimal in the order of magnitude in $N$ but until recently there was no proof of this conjecture (although it was some times quoted as a proven fact). In the meantime there exist explicit constructions of infinite sequences with optimal order of $L_p$ discrepancy in the sense of the general lower bound \eqref{lbdlpdiseq}. These constructions will be presented in Section~\ref{secHOdS}.

\section{Discrepancy of digital sequences}

In the following we give the definition of digital sequences in prime bases $b$. For the general definition we refer to \cite[Section~4.3]{niesiam}. From now on let $b$ be a prime number and let $\FF_b$ be the finite field of order $b$. We identify $\FF_b$ with the set of integers $\{0,1,\ldots,b-1\}$ equipped with the usual arithmetic operations modulo $b$.

\begin{definition}[Niederreiter 1987]
A digital sequence is constructed in the following way:
\begin{itemize}
\item choose $s$ infinite matrices $C_1,\ldots, C_s \in \FF_b^{\NN \times \NN}$;
\item for $n \in \NN_0$ of the form $n = n_0 + n_1 b + n_2 b^2+ \cdots$ and $j=1,2,\ldots,s$ compute (over $\FF_b$) the matrix-vector products 
\begin{equation*}
C_j \left(
\begin{array}{l}
n_0\\
n_1\\
n_2\\
\vdots
\end{array}\right) =:\left(
\begin{array}{l}
x_{n,j,1}\\
x_{n,j,2}\\
x_{n,j,3}\\
\vdots
\end{array}
\right);
\end{equation*}
\item put
\begin{equation*}
x_{n,j} = \frac{x_{n,j,1}}{b} + \frac{x_{n,j,2}}{b^2} + \frac{x_{n,j,3}}{b^3}+\cdots \ \ \ \mbox{ and }\ \ \   \boldsymbol{x}_n = (x_{n,1}, \ldots, x_{n,s}).
\end{equation*}
\end{itemize}
The resulting sequence $\cS(C_1,\ldots,C_s)=(\bsx_n)_{n \ge 0}$ is called a {\it digital sequence over $\FF_b$} and $C_1,\ldots,C_s$ are called the {\it generating matrices} of the digital sequence.
\end{definition}

\subsection{A metrical result}

It is known that almost all digital sequences in a fixed dimension $s$ are low-discrepancy sequences, up to some $\log\log N$-term. The ``almost all'' statement is with respect to a natural probability measure on the set of all $s$-tuples $(C_1,\ldots,C_s)$ of $\NN \times \NN$ matrices over $\FF_b$. For the definition of this probability measure we refer to \cite[p.~107]{LP14}. 

\begin{theorem}[Larcher 1998, Larcher \& Pillichshammer 2014]\label{thmmetric}
Let $\varepsilon>0$. For almost all $s$-tuples $(C_1,\ldots,C_s)$ with $C_j \in \FF_b^{\NN \times \NN}$ the corresponding digital sequences $\cS=\cS(C_1,\ldots,C_s)$ satisfy 
$$D_N^*(\cS) \lesssim_{b,s,\varepsilon} \frac{(\log N)^s (\log \log N)^{2+\varepsilon}}{N} \ \ \ \forall N \ge 2$$ and  
$$D_N^*(\cS) \ge c_{b,s} \frac{(\log N)^s \log \log N}{N} \ \ \ \mbox{infinitely often.}$$
\end{theorem}

The upper estimate has been shown by Larcher in \cite{L98} and a proof for the lower bound can be found in \cite{LP14}.

A corresponding result for the sub-class of so-called digital Kronecker sequences can be found in \cite{L95} (upper bound) and \cite{LP14a} (lower bound). These results correspond to metrical discrepancy bounds for classical Kronecker sequences by Beck~\cite{be}.

The question now arises whether there are $s$ tuples $(C_1,\ldots,C_s)$ of generating matrices such that the resulting digital sequences are low-discrepancy sequences and, if the answer is {\it yes}, which properties of the matrices guarantee low discrepancy. Niederreiter found out that this depends on a certain linear independence structure of the rows of the matrices $C_1,\ldots,C_s$. This leads to the concept of digital $(t,s)$-sequences.

\subsection{Digital $(t,s)$-sequences}\label{digtssequ}

For $C \in \FF_b^{\NN \times \NN}$ and $m \in \NN$ denote by $C(m)$ the left upper $m \times m$ submatrix of $C$.

For technical reasons one often assumes that the generating matrices $C_1,\ldots,C_s$ satisfy the following condition: let $C_j = (c_{k,\ell}^{(j)})_{k, \ell \in \mathbb{N}}$, then for each $\ell \in \mathbb{N}$ there exists a $K(\ell) \in \mathbb{N}$ such that $c_{k,\ell}^{(j)} = 0$ for all $k > K(\ell)$. This condition, which is condition (S6) in \cite[p.72]{niesiam}, guarantees that the components of the elements of a digital sequence have a finite digit expansion in base $b$. For the rest of the paper we tacitly assume that this condition is satisfied. (We remark that in order to include new important constructions to the concept of digital $(t,s)$-sequences, Niederreiter and Xing~\cite{NX1996,NX_book} use a truncation operator to overcome the above-mentioned technicalities. Such sequences are sometimes called $(t,s)$-sequences {\it in the broad sense}.)

\begin{definition}[Niederreiter]
Given $C_1,\ldots,C_s \in \FF_b^{\NN \times \NN}$. If there exists a number $t \in \NN_0$ such that for every $m \ge t$ and for all $d_1,\ldots,d_s\ge 0$ with $d_1+\cdots+d_s=m-t$ the 
$$\left.\begin{array}{l}
 \mbox{first $d_1$ rows of $C_1(m)$,}\\
 \mbox{first $d_2$ rows of $C_2(m)$,}\\
 \ldots \\
 \mbox{first $d_s$ rows of $C_s(m)$,}
\end{array}
\right\} \ \ \mbox{are linearly independent over $\FF_b$,}
$$
then the corresponding digital sequence $\cS(C_1,\ldots,C_s)$ is called a {\it digital $(t,s)$-sequence over $\FF_b$}.
\end{definition}

The technical condition from the above definition guarantees that every $b^m$-element sub-block $(\bsx_{k b^m},\bsx_{k b^m+1},\ldots,\bsx_{(k+1) b^m-1})=:\cS_{m,k}$ of the digital sequence, where $m \ge t$ and $k \in \NN_0$, is a $(t,m,s)$-net in base $b$, i.e., every so-called elementary $b$-adic interval of the form $$J=\prod_{j=1}^s \left[\frac{a_j}{b^{d_j}}, \frac{a_j+1}{b^{d_j}}\right) \ \ \ \ \ \ \mbox{with}\  {\rm Volume}(J)=b^{t-m}$$ contains the right share of elements from $\cS_{m,k}$, which is exactly $b^t$. For more information we refer to \cite[Chapter~4]{DP2010} and \cite[Chapter~4]{niesiam}. Examples for digital $(t,s)$-sequences are generalized Niederreiter sequences which comprise the concepts of Sobol'-, Faure- and original Niederreiter-sequences, Niederreiter-Xing sequences, \ldots. We refer to \cite[Chapter~8]{DP2010} for a collection of constructions and for further references. An overview of the constructions of Niederreiter and Xing can also be found in \cite[Chapter~8]{NX_book}.

It has been shown by Niederreiter~\cite{N87} that every digital $(t,s)$-sequence is a low-discrepancy sequence. The following result holds true:

\begin{theorem}[Niederreiter 1987]\label{thmNie87}
For every digital $(t,s)$-sequence $\cS$ over $\FF_b$ we have $$D_N^*(\cS) \le c_{s,b} \, b^t \, \frac{(\log N)^s}{N}+ O\left(\frac{(\log N)^{s-1}}{N}\right).$$
\end{theorem}

Later several authors worked on improvements of the implied quantity $c_{s,b}$, e.g. \cite{FK,K06}. The currently smallest values for $c_{s,b}$ were provided by Faure and Kritzer~\cite{FK}. More explicit versions of the estimate in Theorem~\ref{thmNie87} can be found in \cite{FL12,FL14,FL15}. For a summary of these results one can also consult \cite[Section~4.3]{FKP}. 

\begin{remark}
Theorem~\ref{thmNie87} in combination with the lower bound in Theorem~\ref{thmmetric} shows that the set of $s$-tuples $(C_1,\ldots,C_s)$ of matrices that generate a digital $(t,s)$-sequence is a set of measure zero. 
\end{remark}

Remember that the exact order of optimal star discrepancy of infinite sequences is still unknown (except for the one-dimensional case). From this point of view it might be still possible that Niederreiter's star discrepancy bound in Theorem~\ref{thmNie87} could be improved in the order of magnitude in $N$. However, it has been shown recently by Levin~\cite{Lev} that this is not possible in general. In his proofs Levin requires the concept of $d$-admissibility. He calls a sequence $(\bsx_n)_{n \ge 0}$ in $[0,1)^s$ {\it $d$-admissible} if $$\inf_{n>k \ge0} \|n\ominus k\|_b \|\bsx_n \ominus \bsx_k\|_b \ge b^{-d},$$ where $\log_b \|x\|_b =\lfloor \log_b x\rfloor$ and $\ominus$ is the $b$-adic difference. Roughly speaking, this means that the $b$-adic distance between elements from the sequence whose indices are close is not too small.

\begin{theorem}[Levin 2017]
Let $\cS$ be a $d$-admissible $(t,s)$-sequence. Then 
$$D_N^{\ast}(\cS) \ge c_{s,t,d} \frac{(\log N)^s}{N}\ \ \ \ \mbox{infinitely often.}$$
\end{theorem}

In his paper, Levin gave a whole list of digital $(t,s)$-sequences that have the property of being $d$-admissible for certain $d$. This list comprises the concepts of generalized Niederreiter sequences (which includes Sobol'-, Faure- and original Niederreiter-sequences), Niederreiter-Xing sequences, \ldots. For a survey of Levin's result we also refer to \cite{KaSt}. It should also be mentioned that there is one single result by Faure~\cite{fau95} from the year 1995 who already gave a lower bound for a particular digital $(0,2)$-sequence (in dimension 2) which is also of order $(\log N)^2/N$.

Levin's results \cite{lev0,Lev} are important contributions to the grand problem in discrepancy theory (cf. Conjecture~\ref{con1}). But they only cover the important sub-class of admissible $(t,s)$-sequences and allow no conclusion for arbitrary (including non-digital) sequences.

\subsection{Digital $(0,1)$-sequences over $\FF_2$}

\begin{figure}[ht]
\begin{center}
\includegraphics[width=10cm]{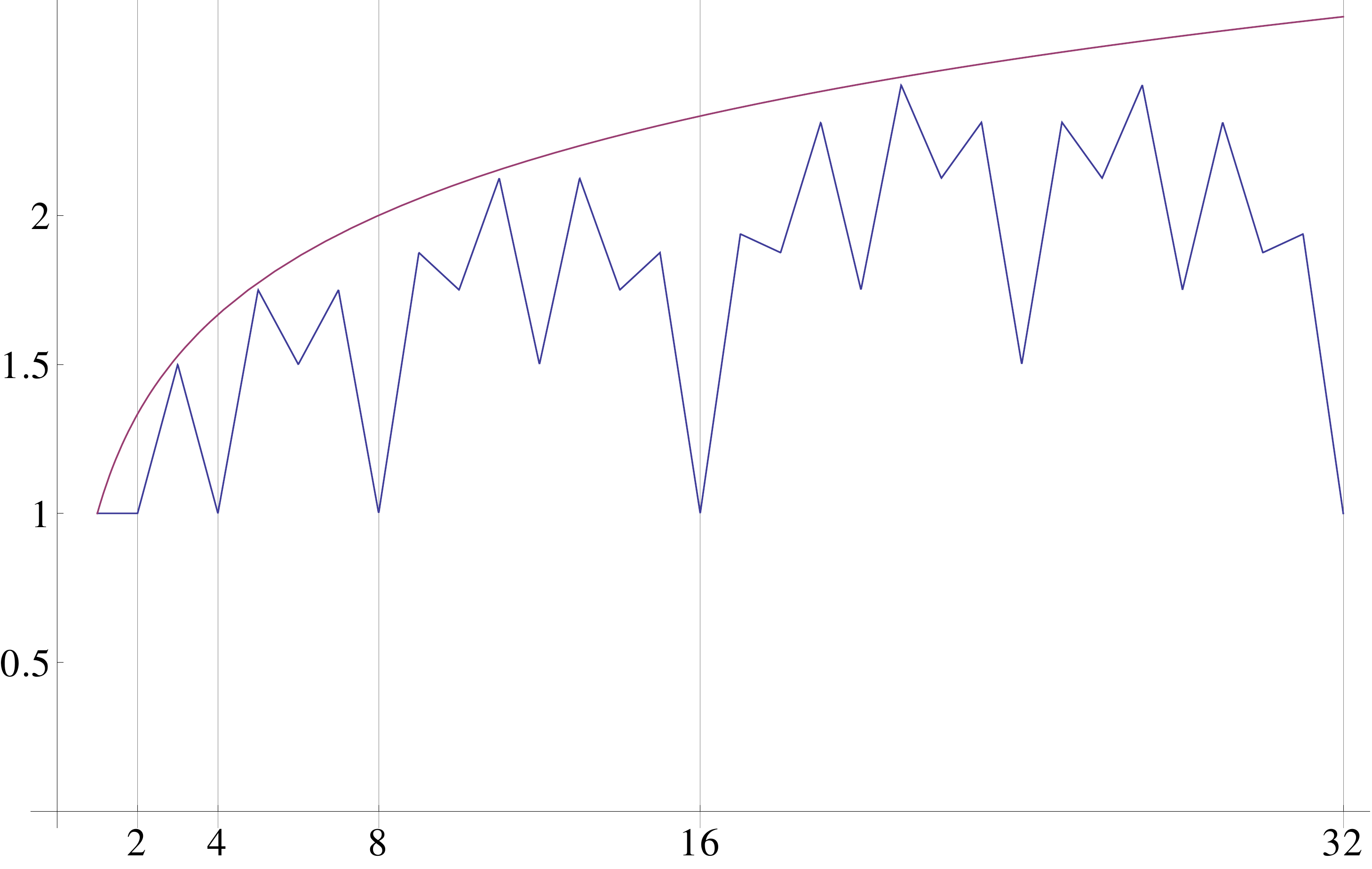}
\end{center}
\caption{$N D_N^{\ast}(\cS(I))$ compared with $\frac{\log N}{3 \log 2} +1$ (red line) for $N=2,3,\ldots,32$; if $N$ is a power of two, then $N D_N^{\ast}(\cS(I))=1$}\label{fig1}
\end{figure}

In this sub-section we say a few words about the discrepancy of digital $(0,1)$-sequence over $\FF_2$, because in this case exact results are known. Let $b=2$ and let $I$ be the $\NN \times \NN$ identity matrix, that is, the matrix whose entries are 0 except for the entries on the main-diagonal which are 1. The corresponding one-dimensional digital sequence $\cS(I)$ is the van der Corput sequence in base $2$ and in fact, it is also a digital $(0,1)$-sequence over $\FF_2$. The following is known: among all digital $(0,1)$-sequences over $\FF_2$ the van der Corput sequence, which is the prototype of all digital constructions and whose star discrepancy is very well studied, has the worst star discrepancy; see \cite[Theorem~2]{P04}. More concretely, for every $\NN \times \NN$ matrix $C$ which generates a digital $(0,1)$-sequence $\cS(C)$ over $\FF_2$ we have   
\begin{equation}\label{estvdcdisc}
D_N^*(\cS(C)) \le D_N^*(\cS(I)) \le \left\{
\begin{array}{l}
\left(\frac{\log N}{3 \log 2} +1\right)\frac{1}{N},\\[0.7em]
\frac{S_2(N)}{N},  
\end{array}\right.
\end{equation}
where $S_2(N)$ denotes the {\it dyadic sum-of-digits function} of the integer $N$. The first bound on $D_N^*(\cS(I))$ is a result of B\'{e}jian and Faure~\cite{befa77}. The factor $1/(3\log 2)$ conjoined with the $\log N$-term is known to be best possible, in fact, $$\limsup_{N \rightarrow \infty} \frac{N D_N^{\ast}(\cS(I))}{\log N}= \frac{1}{3 \log 2}.$$ 
(The corresponding result for van der Corput sequences in arbitrary base can be found in \cite{Fau1981,fau07,K05}.) 
However, also the second estimate in terms of the dyadic sum-of-digits function, which follows easily from the proof of \cite[Theorem~3.5 on p. 127]{kuinie}, is very interesting. It shows that the star discrepancy of the van der Corput sequence (and of any digital $(0,1)$-sequence) is not always close to the high level of order $\log N/N$. If $N$ has only very few dyadic digits different from zero, then the star discrepancy is very small. For example, if $N$ is a power of two, then $S_2(N)=1$ and therefore $D_N^{\ast}(\cS(I))\le 1/N$. The bound in \eqref{estvdcdisc} is demonstrated in Fig.~\ref{fig1}.

While the star discrepancy of any digital $(0,1)$-sequence over $\FF_2$ is of optimal order with respect to \eqref{lpseqschm} this fact is not true in general for the $L_p$ discrepancies with finite parameter $p$. For example, for the van der Corput sequence we have for all $p \in [1,\infty)$ $$\limsup_{N \rightarrow \infty} \frac{N L_{p,N}(\cS(I))}{\log N}= \frac{1}{6 \log 2},$$ see \cite{P04}. Hence the $L_p$ discrepancy of the van der Corput sequence is at least of order of magnitude $\log N/N$ for infinitely many $N$. Another example, to be found in \cite{DLP}, is the digital $(0,1)$-sequence generated by the matrix $$U=\left(\begin{array}{cccc} 
               1 & 1 & 1 & \ldots\\
               0 & 1 & 1 & \ldots\\
               0 & 0 & 1 & \ldots\\
              \multicolumn{4}{l}\dotfill 
              \end{array}
\right)$$ for which we have, with some positive real $c>0$, $$\limsup_{N \rightarrow \infty} \frac{N L_{2,N}(\cS(U))}{\log N} \ge c >0.$$ 

More information on the discrepancy of digital $(0,1)$-sequences can be found in the survey articles \cite{FK14,FKP} and the references therein.

The results in dimension one show that, in general, the $L_p$ discrepancy of digital sequences does not match the general lower bound \eqref{lbdlpdiseq} from Roth-Schmidt-Pro{\u\i}nov. Hence, in order to achieve the assumed optimal order of magnitude $(\log N)^{s/2}/N$ for the $L_p$ discrepancy with digital sequences, if at all possible, one needs more demanding properties on the generating matrices. This leads to the concept of higher order digital sequences.

\subsection{Digital sequences with optimal order of $L_p$ discrepancy}\label{secHOdS}

So-called {\it higher order digital sequences} have been introduced by Dick \cite{D07,D08} in 2007 with the aim to achieve optimal convergence rates for QMC rules applied to sufficiently smooth functions. For the definition of higher order digital sequences and for further information and references we refer to \cite[Chapter~15]{DP2010} or to \cite{DKS}.

For our purposes it suffices to consider higher order digital sequences of order 2. We just show how such sequences can be constructed: to this end let $d:=2 s$ and let $C_1, \ldots, C_{d} \in \FF_2^{\NN \times \NN}$ be generating matrices of a digital $(t,d)$-sequence in dimension $d$, for example a generalized Niederreiter sequence. Let $\vec{c}_{j,k}$ denote the $k^{{\rm th}}$ row-vector of the matrix $C_j$. Now define $s$ matrices $E_1,\ldots, E_s$ in the following way: the row-vectors of $E_j$ are given by 
\begin{equation*}
\vec{e}_{j,2 u + v} = \vec{c}_{2 (j-1) + v, u+1} \ \ \ \mbox{for $j \in \{1,2,\ldots,s\}$, $u \in \N_0$ and $v \in \{1,2\}$.}
\end{equation*}
We illustrate the construction for $s=1$. Then $d=2$ and
$$C_1=\left(\begin{array}{c} 
             \vec{c}_{1,1}\\
             \vec{c}_{1,2}\\
             \vdots 
            \end{array}\right), \ C_2=\left(\begin{array}{c} 
             \vec{c}_{2,1}\\
             \vec{c}_{2,2}\\
             \vdots 
            \end{array}\right) \ \Rightarrow \ 
           E_1=\left(\begin{array}{c} 
             \vec{c}_{1,1}\\
             \vec{c}_{2,1}\\
             \vec{c}_{1,2}\\
             \vec{c}_{2,2}\\
             \vdots 
            \end{array}\right).$$
This procedure is called \emph{interlacing} (here the so-called ``interlacing factor'' is $2$).

The following theorem has been shown in \cite{DHMP17a}.

\begin{theorem}[Dick, Hinrichs, Markhasin \& Pillichshammer 2017]\label{thm_main2}
Assume that $E_1,\ldots,E_s \in \FF_2^{\NN \times \NN}$ are constructed with the interlacing principle as given above. Then for the corresponding digital sequence $\cS=\cS(E_1,\ldots,E_s)$ we have  
$$
L_{p,N}(\cS) \lesssim_{p,s} 2^{2t} \frac{\left(\log N\right)^{s/2}}{N}\ \ \ \ \mbox{ for all $N\ge 2$ and all $1 \le p < \infty$.}
$$
\end{theorem}

This theorem shows, in a constructive way, that the lower bound \eqref{lbdlpdiseq} from Roth-Schmidt-Pro{\u\i}nov is best possible in the order of magnitude in $N$ for all parameters $p\in (1,\infty)$. Furthermore, the constructed digital sequences have optimal order of $L_p$ discrepancy simultaneously for all $p \in (1,\infty)$.

For $p=2$ there is an interesting improvement, although this improvement requires higher order digital sequences of order 5 (instead of order 2). For such sequences $\cS$ it has been shown in \cite{DP14b} that $$L_{2,N}(\cS) \lesssim_s \frac{(\log N)^{(s-1)/2}}{N} \, \sqrt{S_2(N)}\ \ \ \ \mbox{ for all $N\ge 2$.}$$ The dyadic sum-of-digit function of $N$ is in the worst-case of order $\log N$ and then the above $L_2$ discrepancy bound is of order of magnitude $(\log N)^{s/2}/N$. But if $N$ has very few non-zero dyadic digits, for example if it is a power of 2, then the bound on the $L_2$ discrepancy becomes $(\log N)^{(s-1)/2}/N$ only.

The proof of Theorem~\ref{thm_main2} uses methods from harmonic analysis, in particular the estimate of the $L_p$ norm of the discrepancy function is based on the following  Littlewood-Paley type inequality: for $p \in (1,\infty)$ and $f\in L_p([0,1]^s)$ we have
\begin{equation}\label{LiPal} 
\|f\|_{L_p([0,1]^s)} \lesssim_{p,s} \sum_{\bsj\in\NN_{-1}^s} 2^{2|\bsj|(1-1/\bar{p})}\left(\sum_{\bsm\in\mathbb{D}_{\bsj}} |\langle f,h_{\bsj,\bsm}\rangle|^{\bar{p}}\right)^{2/\bar{p}}, 
\end{equation}
where $\bar{p} = \max(p,2)$, $\NN_{-1}=\NN \cup \{-1,0\}$, for $\bsj=(j_1,\ldots,j_s)$, $\mathbb{D}_{\bsj}=\mathbb{D}_{j_1}\times \ldots \times \mathbb{D}_{j_s}$, where $\mathbb{D}_j=\{0,1,\ldots,2^j-1\}$, $|\bsj|=\max(j_1,0)+\cdots+\max(j_s,0)$, and, for $\bsm \in \mathbb{D}_{\bsj}$,  $h_{\bsj,\bsm}(\bsx)=h_{j_1,m_1}(x_1)\cdots h_{j_s,m_s}(x_s)$, where $h_{j,m}$ is the $m^{{\rm th}}$ dyadic Haar function on level $j$. See \cite{DHMP17a} and \cite{M13b}. The $L_2$ inner products $\langle f,h_{\bsj,\bsm}\rangle$ are the so-called Haar coefficients of $f$. Inequality \eqref{LiPal} is used for the local discrepancy function of digital sequences which then requires tight estimates of the Haar coefficients of the local discrepancy function. For details we refer to \cite{DHMP17a}.

With the same method one can also handle the quasi-norm of the local discrepancy function in Besov spaces and Triebel-Lizorkin spaces with dominating mixed smoothness. One reason why Besov spaces and Triebel-Lizorkin spaces are interesting in this context is that they form natural scales of function spaces including the $L_p$-spaces and Sobolev spaces of dominating mixed smoothness (see, e.g., \cite{T10}). The study of discrepancy in these function spaces has been initiated by Triebel \cite{T10,T10a} in 2010. Further results (for finite sequences) can be found in \cite{H10,M13a,M13b,M13c,M15} and (for infinite sequences in dimension one) in \cite{K15}.  In \cite[Theorem~3.1 and 3.2]{DHMP17b} general lower bounds on the quasi-norm of the local discrepancy function in Besov spaces and Triebel-Lizorkin spaces with dominating mixed smoothness in the sense of the result of Roth-Schmidt-Pro{\u\i}nov in Eq. \eqref{lbdlpdiseq} are shown. Furthermore, these lower bounds are optimal in the order of magnitude in $N$, since matching upper bounds are obtained for infinite order two digital sequences as constructed above. For details we refer to \cite{DHMP17b}.

\subsection{Intermediate norms of the local discrepancy function}

While the quest for the exact order of the optimal $L_p$ discrepancy of infinite sequences in arbitrary dimension is now solved for finite parameters $p\in (1,\infty)$ the situation for the cases $p\in \{1,\infty\}$ remains open. In this situation, Bilyk, Lacey, Parissis and Vagharshakyan \cite{BLPV09} studied the question of what happens in intermediate spaces ``close'' to $L_{\infty}$. Two standard examples of such spaces are:
\begin{itemize}
\item {\it Exponential Orlicz space}: for the exact definition of the corresponding norm $\|\cdot\|_{\exp(L^\beta)}$, $\beta>0$, we refer to \cite{BLPV09,BM15,DHMP17b}. There is an equivalence which shows the relation to the $L_p$ norm, which is stated for any $\beta > 0$, 
\begin{equation*}
\| f \|_{\exp(L^\beta)} \asymp \sup_{p > 1} p^{-\frac{1}{\beta}} \|f \|_{L_p([0,1]^s)}.
\end{equation*}
This equivalence suggests that the study of discrepancy with respect to the exponential Orlicz norm is related to the study of the dependence of the constant appearing in the $L_p$ discrepancy bounds on the parameter $p$. The latter problem is also studied in \cite{skr12}. 
\item {\it BMO space} (where BMO stands for ``bounded mean oscillation''): the definition of the corresponding semi-norm uses Haar functions and is given as  $$\|f\|_{{\rm BMO}^s}^2 =\sup_{U \subseteq [0,1)^s} \frac{1}{\lambda_s(U)} \sum_{\bsj \in \mathbb{N}_0^s} 2^{|\bsj|} \sum_{\bsm \in \mathbb{D}_{\bsj}\atop \textsf{supp}(h_{\bsj,\bsm}) \subseteq U}|\langle f, h_{\bsj,\bsm}\rangle|^2,$$ where the supremum is extended over all measurable subsets $U$ from $[0,1)^s$. See again \cite{BLPV09,BM15,DHMP17b} and the references therein for more information.
\end{itemize}

Exponential Orlicz norm and BMO semi-norm of the local discrepancy function for finite point sets have been studied in \cite{BLPV09} (in dimension $s=2$) and in \cite{BM15} (in the general multi-variate case). For infinite sequences we have the following results which have been shown in \cite{DHMP17b}:

\begin{theorem}[Dick, Hinrichs, Markhasin \& Pillichshammer 2017]\label{thm_main3}
Assume that $E_1,\ldots,E_s \in \FF_2^{\NN \times \NN}$ are constructed with the interlacing principle as given in Section~\ref{secHOdS}. Then for the corresponding digital sequence $\cS=\cS(E_1,\ldots,E_s)$ we have  
$$
\|\Delta_{\cS_N}\|_{\exp(L^{\beta})} \lesssim_s \frac{(\log N)^{s-\frac{1}{\beta} }}{N} \ \ \ \mbox{ for all $N \ge 2$ and for all $\frac{2}{s-1} \le \beta < \infty$}
$$
and 
\begin{equation}\label{bmo:ubd}
\|\Delta_{\cS_N}\|_{{\rm BMO}^s} \lesssim_s \frac{(\log N)^{\frac{s}{2}}}{N} \ \ \ \mbox{ for all $N \ge 2$.}
\end{equation} 
\end{theorem}

A matching lower bound in the case of exponential Orlicz norm on $\|\Delta_{\cS_N}\|_{\exp(L^{\beta})}$ in arbitrary dimension is currently not available and seems to be a very difficult problem, even for finite sequences (see \cite[Remark after Theorem~1.3]{BM15}; for matching lower and upper bounds for finite sequences in dimension $s=2$ we refer to \cite{BLPV09}). On the other hand, the result \eqref{bmo:ubd} for the BMO semi-norm is best possible in the order of magnitude in $N$. A general lower bound in the sense of Roth-Schmidt-Pro{\u\i}nov's result \eqref{lbdlpdiseq} for the $L_p$ discrepancy has been shown in \cite[Theorem~2.1]{DHMP17b} and states that for every $s \in \NN$ there exists a $c_s>0$ such that for every infinite sequence $\cS$ in $[0,1)^s$ we have 
\begin{equation}\label{bmo:lowbd}
\|\Delta_{\cS_N}\|_{{\rm BMO}^s} \ge c_s  \frac{(\log N)^{\frac{s}{2}}}{N}\ \ \ \mbox{ infinitely often.}
\end{equation}

\section{Discussion of the asymptotic discrepancy estimates}

We restrict the following discussion to the case of star discrepancy. We have seen that the star discrepancy of digital sequences, and therefore QMC rules which are based on digital sequences, can achieve error bounds of order of magnitude $(\log N)^s/N$. At first sight this seems to be an excellent result. However, the crux of these, in an asymptotic sense, optimal results, lies in the dependence on the dimension $s$. If we consider the function $x \mapsto  (\log x)^s/x$, then one can observe, that this function is increasing up to $x={\rm e}^s$ and only then it starts to decrease to 0 with the asymptotic order of almost $1/x$. This means, in order to have meaningful error bounds for QMC rules one requires finite sequences with at least ${\rm e}^s$ many elements or even larger. But ${\rm e}^s$ is already huge, even for moderate dimensions $s$. For example, if $s=200$, then ${\rm e}^s \approx 7.2 \times 10^{86}$ which exceeds the estimated number of atoms in our universe (which is $\approx10^{78}$).

As it appears, according to the classical theory with its excellent asymptotic results, QMC rules cannot be expected to work for high-dimensional functions. However, there is numerical evidence, that QMC rules can also be used in these cases. The work of Paskov and Traub~\cite{PT95} from 1995 attracted much attention in this context. They considered a real world problem from mathematical finance which resulted in the evaluation of several 360 dimensional integrals and reported on their successful use of Sobol' and Halton-sequences in order to evaluate these integrals.  

Of course, it is now the aim of theory to the explain, {\it why} QMC rules also work for high-dimensional problems. One stream of research is to take the viewpoint of {\it Information Based Complexity (IBC)} in which also the dependence of the error bounds (discrepancy in our case) on the dimension $s$ is studied. A first remarkable, and at that time very surprising result, has been established by Heinrich, Novak, Wasilkowski and Wo\'{z}niakowski~\cite{HNWW} in 2001. 

\begin{theorem}[Heinrich, Novak, Wasilkowski \& Wo\'{z}niakowski 2001]\label{thmHNWW}
For all $N,s \in\NN$ there exist finite sequences $\cS_N$ of $N$ elements in $[0,1)^s$ such that $$D_N^{\ast}(\cS_N) \lesssim \sqrt{\frac{s}{N}},$$ where the implied constant is absolute, i.e., does neither depend on $s$, nor on $N$. 
\end{theorem}

In 2007 Dick~\cite{D07a} extended this result to infinite sequences (in infinite dimension).
   
In IBC the information complexity is studied rather then direct error bounds. In the case of star discrepancy the {\it information complexity}, which is then also called the {\it inverse of star discrepancy}, is, for some error demand $\varepsilon \in (0,1]$ and dimension $s$, given as $$N^{\ast}(\varepsilon,s)=\min\{N \in \NN \ : \ \exists \ \cS_N \subseteq [0,1)^s \ \mbox{with $|\cS_N|=N$ and $D_N^{\ast}(\cS_N) \le \varepsilon$}\}.$$
From Theorem~\ref{thmHNWW} one can deduce that $$N^{\ast}(\varepsilon,s) \lesssim s \varepsilon^{-2}$$ and this property is called {\it polynomial tractability} with $\varepsilon$-exponent $2$ and $s$-exponent 1. In 2004 Hinrichs~\cite{hin} proved that there exists a positive $c$ such that $N^{\ast}(\varepsilon,s) \ge c s \varepsilon^{-1}$ for all $s$ and all small enough $\varepsilon>0$. Combining these results we see, that {\it the inverse of the star discrepancy depends (exactly) linearly on the dimension $s$} (which is the programmatic title of the paper \cite{HNWW}). The exact dependence of the inverse of the star discrepancy on $\varepsilon^{-1}$ is still unknown and seems to be a very difficult problem. In 2011 Aistleitner~\cite{Aist} gave a new proof of the result in  Theorem~\ref{thmHNWW} from which one can obtain an explicit constant in the star discrepancy estimate. He proved that there exist finite sequences $\cS_N$ of $N$ elements in $[0,1)^s$ such that $D_N^{\ast}(\cS_N) \le 10 \sqrt{s/N}$ and hence $N^{\ast}(\varepsilon,s)\le 100 s \varepsilon^{-2}$. Recently Gnewuch and Hebbinghaus (private communication) improved these implied constants to $D_N^{\ast}(\cS_N) \le (2.5287\ldots) \times \sqrt{s/N}$  and hence $N^{\ast}(\varepsilon,s)\le  (6.3943\ldots)\times s \varepsilon^{-2}$. 

For a comprehensive introduction to IBC and tractability theory we refer to the three volumes \cite{NW08,NW10,NW12} by Novak and Wo\'{z}niakowski.

Unfortunately, the result in Theorem~\ref{thmHNWW} is a pure existence result and until now no concrete point set is known whose star discrepancy satisfies the given upper bound. Motivated by the excellent asymptotic behavior it may be obvious to consider digital sequences also in the context of tractability. This assumption is supported by a recent metrical result for a certain subsequence of a digital Kronecker sequence. In order to explain this result we need some notation:
\begin{itemize}
\item Let $\FF_b((t^{-1}))$ be the field of {\it formal Laurent series} over $\FF_b$ in the variable $t$:
$$\FF_b((t^{-1})) = \left\{ \sum_{i=w}^{\infty} g_i \, t^{-i} ~ : ~ w \in \ZZ, \forall i: g_i \in \FF_b \right\}.$$
\item For $g \in \FF_b((t^{-1}))$ of the form $g = \sum_{i=w}^{\infty}g_i \, t^{-i}$ define the ``fractional part'' $$\{g\}:=\sum_{i=\max\{w,1\}}^{\infty} g_i \, t^{-i}.$$ 
\item Every $n \in \NN_0$ with $b$-adic expansion $n=n_0+n_1 b+\cdots + n_r b^{r}$, where $n_i \in \{0,\ldots,b-1\}$, is associated in the natural way with the polynomial
$$n \cong n_0 + n_1t + \cdots + n_rt^r \in \FF_b[t].$$
\end{itemize}

Now a digital Kronecker sequence is defined as follows:

\begin{definition}
Let $\bsf=(f_1,\ldots,f_s) \in \FF_b((t^{-1}))^s$. Then the sequence $\mathcal{S}(\bsf)=(\bsy_n)_{n\geq 0}$ given by 
\begin{align*}
\bsy_n:=\{n\bsf\}_{|t=b}=(\{nf_1\}_{|t=b},\ldots,\{nf_s\}_{|t=b})
\end{align*}
is called a {\it digital Kronecker sequence over $\FF_b$}.
\end{definition}

It can be shown that digital Kronecker sequences are examples of digital sequences where the generating matrices are Hankel matrices (i.e., constant ascending skew-diagonals) whose entries are the coefficients of the Laurent series expansions of $f_1,\ldots,f_s$; see, e.g., \cite{LN93,niesiam}. Neum\"uller and Pillichshammer~\cite{NP18} studied a subsequence of digital Kronecker sequences. For $\bsf \in \FF_b((t^{-1}))^s$ consider 
$\widetilde{\mathcal{S}}(\bsf)=(\bsy_n)_{n \ge 0}$ where $$\bsy_n=\{t^{n}\bsf\}_{|t=b}=(\{t^{n}f_1\}_{|t=b},\ldots,\{t^{n}f_s\}_{|t=b}).$$
With a certain natural probability measure on $\FF_b((t^{-1}))^s$ the following metrical result can be shown:

\begin{theorem}[Neum\"uller \& Pillichshammer 2018]\label{thmNP18}
Let $s \ge 2$. For every $\delta \in (0,1)$ we have
\begin{align}\label{estNP}
D^*_N(\widetilde{\cS}(\bsf)) \lesssim_{b,\delta} \sqrt{\frac{s\log s}{N}}\, \log N \ \ \ \mbox{ for all $N \ge 2$}
\end{align}
with probability at least $1-\delta$, where the implied constant $C_{b,\delta} \asymp_b \log \delta^{-1}$.
\end{theorem} 
  
The estimate \eqref{estNP} is only slightly weaker than the bound in Theorem~\ref{thmHNWW}. The additional $\log N$-term comes from the consideration of infinite sequences. Note that the result holds for all $N \ge 2$ simultaneously. One gets rid of this $\log N$-term when one considers only finite sequences as in Theorem~\ref{thmHNWW}; see \cite[Theorem~3]{NP18}. Furthermore, we remark that Theorem~\ref{thmNP18} corresponds to a result for classical Kronecker sequences which has been proved by L\"obbe~\cite{loeb}.  

\section{Weighted discrepancy of digital sequences}

Another way to explain the success of QMC rules for high-dimensional problems is the study of so-called weighted function classes. This study, initiated by Sloan and Wo\'{z}niakowski~\cite{SW98} in 1998, is based on the assumption that functions depend differently on different variables and groups of variables when the dimension $s$ is large. This different dependence should be reflected in the error analysis. For this purpose Sloan and Wo\'{z}niakowski proposed the introduction of weights that model the dependence of the functions on different coordinate directions. In the context of discrepancy theory this led to the introduction of weighted $L_p$ discrepancy. Here we restrict ourselves to the case of weighted star discrepancy:

In the following let $\bsgamma=(\gamma_1,\gamma_2,\gamma_3,\ldots)$ be a sequence of positive reals, the so-called weights. Let $[s]:=\{1,2,\ldots,s\}$ and for ${\mathfrak u} \subseteq [s]$ put $$\gamma_{{\mathfrak u}}:=\prod_{j \in {\mathfrak u}} \gamma_j.$$

\begin{definition}[Sloan \& Wo\'{z}niakowski 1998]
For a sequence $\cS$ in $[0,1)^s$ the {\em $\bsgamma$-weighted star discrepancy} is defined as
$$D_{N,{\bsgamma}}^*(\cS):=\sup_{\bsalpha\in [0,1]^s}
\max_{\emptyset\ne {\mathfrak u}\subseteq [s]} \gamma_{\mathfrak
  u}|\Delta_{\cS_N}(\bsalpha_{\mathfrak u},\bsone)|,$$
where for $\bsalpha=(\alpha_1,\ldots,\alpha_s) \in [0,1]^s$ and for
${\mathfrak u} \subseteq [s]$ we put $(\bsalpha_{{\mathfrak u}},\bsone)=(y_1,\ldots,y_s)$
with $y_j=\alpha_j$ if $j \in {\mathfrak u}$ and $y_j=1$ if $j \not\in {\mathfrak u}$. 
\end{definition}

\begin{remark}
If $\gamma_j=1$ for all $j \ge 1$, then $D_{N,{\bsgamma}}^*(\cS) = D_N^*(\cS).$
\end{remark}

The relation between weighted discrepancy and error bounds for QMC rules is expressed by means of a {\it weighted Koksma-Hlawka inequality} as follows:
Let $\mathcal{W}_1^{(1,1,\ldots,1)}([0,1]^s)$ be the Sobolev space of functions defined on $[0,1]^s$ that are once differentiable in each variable, and whose derivatives have finite $L_1$ norm. Consider $$\mathcal{F}_{s,1,\bsgamma}=\{f \in \mathcal{W}_1^{(1,1,\ldots,1)}([0,1]^s) \ : \ \|f\|_{s,1,\bsgamma}< \infty\},$$ where $$\|f\|_{s,1,\bsgamma} =|f(\bsone)| + \sum_{\emptyset \not={\mathfrak u} \subseteq [s]} \frac{1}{\gamma_{\mathfrak u}}\left\|\frac{\partial^{|{\mathfrak u}|}}{\partial \bsx_{{\mathfrak u}}}f(\bsx_{{\mathfrak u}},\bsone)\right\|_{L_1}.$$ The $\bsgamma$-weighted star discrepancy of a finite sequence is then exactly the worst-case error of a QMC rule in $\mathcal{F}_{s,1,\bsgamma}$ that is based on this sequence, see \cite{SW98} or  \cite[p.~65]{NW10}. More precisely, we have 
$$\sup_{\|f\|_{s,1,\bsgamma} \le 1} \left|\int_{[0,1]^s} f(\bsx) \rd \bsx - \frac{1}{N} \sum_{\bsx \in \cS_N} f(\bsx)\right|=D_{N,{\bsgamma}}^*(\cS).$$

In IBC again the {\it inverse of weighted star discrepancy}  
$$N_{\bsgamma}^{\ast}(\varepsilon,s) := \min\{N \ : \ \exists \, \cS_N \subseteq [0,1)^s \ \mbox{with $|\cS_N|=N$ and $D_{N,\bsgamma}^*(\cS_N) \le \varepsilon$}\}$$ is studied. The weighted star discrepancy is said to be {\it strongly polynomially tractable} (SPT), if there exist non-negative real numbers $C$ and $\beta$ such that
\begin{equation}\label{defspt}
N_{\bsgamma}^{\ast}(\varepsilon,s) \le C \varepsilon^{-\beta}\ \ \ \mbox{ for all $s\in \NN$ and for all $\varepsilon \in (0,1)$.}
\end{equation}
The infimum $\beta^{\ast}$ over all $\beta > 0$ such that \eqref{defspt} holds is called the $\varepsilon$-exponent of strong polynomial tractability. It should be mentioned, that there are several other notions of tractability which are considered in literature. Examples are polynomial tractability, weak tractability, etc. For an overview we refer to \cite{NW08,NW10,NW12}.

In \cite{HPT18} Hinrichs, Tezuka and the author studied tractability properties of the weighted star discrepancy of several digital sequences.

\begin{theorem}[Hinrichs, Pillichshammer \& Tezuka 2018]\label{thm_HPT}
The weighted star discrepancy of the Halton sequence (where the bases $b_1,\ldots,b_s$ are the first $s$ prime numbers in increasing order) and of Niederreiter sequences achieve SPT with $\varepsilon$-exponent
\begin{itemize}
\item $\beta^{\ast}=1$, which is optimal, if
$$\sum_{j \ge 1} j \gamma_j < \infty, \ \ \ \ \ \ \ \ \ \ \ \mbox{e.g., if $\gamma_j=\frac{1}{j^{2+\delta}}$ with some $\delta>0$;}$$ 
\item $\beta^{\ast} \le 2$, if $$\sup\limits_{s \ge 1} \max\limits_{\emptyset \not= \mathfrak{u} \subseteq [s]} \prod_{j \in \mathfrak{u}} (j \gamma_j) < \infty\ \ \ \ \ \ \ \ \ \ \ \mbox{e.g., if $\gamma_j=\frac{1}{j}$.}$$  
\end{itemize}
\end{theorem}

This result is the currently mildest weight condition for a ``constructive'' proof of SPT of the weighted star discrepancy. Furthermore, it is the first ``constructive'' result which does not require that the weights are summable in order to achieve SPT. By a ``constructive'' result we mean in this context that the corresponding point set can be found or constructed by a polynomial-time algorithm in $s$ and in $\varepsilon^{-1}$. 

To put the result in Theorem~\ref{thm_HPT} into context we recall the currently best ``existence result'' which has been shown by Aistleitner~\cite{Aist2}:

\begin{theorem}[Aistleitner] 
If there exists a $c>0$ such that
\begin{equation*}
\sum_{j=1}^\infty \exp(-c \gamma_j^{-2}) < \infty \ \ \ \ \ \ \ \ \mbox{e.g., if $\gamma_j=\frac{1}{\sqrt{\log j}}$,}
\end{equation*}
then the weighted star discrepancy is SPT with $\varepsilon$-exponent $\beta^* \le 2$. 
\end{theorem}

Obviously the condition on the weights in Aistleitner's ``existence'' result is much weaker then for the ``constructive'' result in Theorem~\ref{thm_HPT}. It is now the task to find sequences whose weighted star discrepancy achieves SPT under the milder weight condition.

\section{Summary}

Digital $(t,s)$-sequences are without doubt the most powerful concept for the construction of low-discrepancy sequences in many settings. Such sequences are very much-needed as sample points for QMC integration rules. They have excellent discrepancy properties in an asymptotic sense when the dimension $s$ is fixed and when $N \rightarrow \infty$:
\begin{itemize}
\item For $p \in [1,\infty)$ there are constructions of digital sequences with $L_p$ discrepancy $$L_p(\cS) \lesssim_{s,p} \frac{(\log N)^{s/2}}{N} \ \ \ \ \mbox{ for all $N \ge 2$ and $p\in [1,\infty)$}$$ and this estimate is best possible in the order of magnitude in $N$ for $p \in (1,\infty)$ according to the general lower bound \eqref{lbdlpdiseq}. 
\item The star discrepancy of digital $(t,s)$-sequences satisfies a bound of the form $$D_N^{\ast}(\cS) \lesssim_s \frac{(\log N)^s}{N} \ \ \ \ \mbox{ for all $N \ge 2$}$$ and this bound is often assumed to be best possible at all.
\item For discrepancy with respect to various other norms digital sequences achieve very good and even optimal results.
\end{itemize}

On the other hand, nowadays one is also very much interested in the dependence of discrepancy on the dimension $s$. This is a very important topic, in particular in order to justify the use of QMC in high dimensions. First results suggest that also in this IBC context digital sequences may perform very well. But here many questions are still open and require further studies. One particularly important question is how sequences can be constructed whose discrepancy achieves some notion of tractability. Maybe digital sequences are good candidates also for this purpose.

\begin{small}
\noindent\textbf{Authors' address:} \\
Friedrich Pillichshammer, Institut f\"{u}r Finanzmathematik und Angewandte Zahlentheorie, Johannes Kepler Universit\"{a}t Linz, Altenbergerstr.~69, 4040 Linz, Austria.\\ \textbf{E-mail:} \texttt{friedrich.pillichshammer@jku.at}
\end{small}
\end{document}